\font\fifteenrm=cmr10 scaled\magstep2 
\font\fifteeni=cmmi10 scaled\magstep2
\font\fifteensy=cmsy10 scaled\magstep2
\font\fifteenbf=cmbx10 scaled\magstep2
\font\fifteentt=cmtt10 scaled\magstep2
\font\fifteenit=cmti10 scaled\magstep2
\font\fifteensl=cmsl10 scaled\magstep2
\font\fifteenam=msam10 scaled\magstep2
\font\fifteenbm=msbm10 scaled\magstep2
\font\fifteenex=cmex10 scaled\magstep2
\font\fifteensc=cmcsc10 scaled\magstep2 
\font\twelverm=cmr10 at 12pt
\font\twelvei=cmmi10 at 12pt
\font\twelvesy=cmsy10 at 12pt
\font\twelvebf=cmbx10 at 12pt
\font\twelvett=cmtt10 at 12pt
\font\twelveit=cmti10 at 12pt
\font\twelvesl=cmsl10 at 12pt
\font\twelveam=msam10 at 12pt
\font\twelvebm=msbm10 at 12pt
\font\twelveex=cmex10 at 12pt
\font\twelvesc=cmcsc10 at 12pt
\font\elevenrm=cmr10 scaled\magstephalf 
\font\eleveni=cmmi10 scaled\magstephalf
\font\elevensy=cmsy10 scaled\magstephalf
\font\elevenbf=cmbx10 scaled\magstephalf
\font\eleventt=cmtt10 scaled\magstephalf
\font\elevenit=cmti10 scaled\magstephalf
\font\elevensl=cmsl10 scaled\magstephalf
\font\elevenam=msam10 scaled\magstephalf
\font\elevenbm=msbm10 scaled\magstephalf
\font\elevenex=cmex10 scaled\magstephalf
\font\elevensc=cmcsc10 scaled\magstephalf
\font\tenrm=cmr10
\font\teni=cmmi10
\font\tensy=cmsy10
\font\tenbf=cmbx10
\font\tentt=cmtt10
\font\tenit=cmti10
\font\tensl=cmsl10
\font\tenam=msam10
\font\tenbm=msbm10
\font\tenex=cmex10
\font\tensc=cmcsc10
\font\ninerm=cmr9
\font\ninei=cmmi9
\font\ninesy=cmsy9
\font\ninebf=cmbx9
\font\ninett=cmtt9
\font\nineit=cmti9
\font\ninesl=cmsl9
\font\nineam=msam9
\font\ninebm=msbm9
\font\nineex=cmex9
\font\ninesc=cmcsc9
\font\eightrm=cmr8
\font\eighti=cmmi8
\font\eightsy=cmsy8
\font\eightbf=cmbx8
\font\eighttt=cmtt8
\font\eightit=cmti8
\font\eightsl=cmsl8
\font\eightam=msam8
\font\eightbm=msbm8
\font\eightex=cmex8
\font\eightsc=cmcsc8
\font\sevenrm=cmr7
\font\seveni=cmmi7
\font\sevensy=cmsy7
\font\sevenbf=cmbx7

\font\sevenam=msam7
\font\sevenbm=msbm7

\font\sixrm=cmr6
\font\sixi=cmmi6
\font\sixsy=cmsy6

\font\sixam=msam6
\font\sixbm=msbm6

\font\fiverm=cmr5
\font\fivei=cmmi5
\font\fivesy=cmsy5

\font\fiveam=msam5
\font\fivebm=msbm5

\font\fourrm=cmr5 at 4pt
\font\fouri=cmmi5 at 4pt
\font\foursy=cmsy5 at 4pt

\font\fouram=msam5 at 4pt
\font\fourbm=msbm5 at 4pt

\skewchar\twelvei='177 \skewchar\eleveni='177\skewchar\teni='177
\skewchar\ninei='177 \skewchar\eighti='177\skewchar\seveni='177 
\skewchar\sixi='177 \skewchar\fivei='177 \skewchar\fouri='177
\skewchar\twelvesy='60 \skewchar\elevensy='60 \skewchar\tensy='60
\skewchar\ninesy='60 \skewchar\eightsy='60 \skewchar\sevensy='60 
\skewchar\sixsy='60 \skewchar\fivesy='60 \skewchar\foursy='60
\newfam\itfam
\newfam\slfam
\newfam\bffam
\newfam\ttfam
\newfam\scfam
\newfam\amfam
\newfam\bmfam
\def\eightbig#1{{\hbox{$\left#1\vbox to 6.5pt{}\voidright $}}}
\def\eightBig#1{{\hbox{$\left#1\vbox to 7.5pt{}\voidright $}}}
\def\eightbigg#1{{\hbox{$\left#1\vbox to 10pt{}\voidright $}}}
\def\eightBigg#1{{\hbox{$\left#1\vbox to 13pt{}\voidright $}}}
\def\ninebig#1{{\hbox{$\left#1\vbox to 7.5pt{}\voidright $}}}
\def\nineBig#1{{\hbox{$\left#1\vbox to 8.5pt{}\voidright $}}}
\def\ninebigg#1{{\hbox{$\left#1\vbox to 11.5pt{}\voidright $}}}
\def\nineBigg#1{{\hbox{$\left#1\vbox to 14.5pt{}\voidright $}}}
\def\tenbig#1{{\hbox{$\left#1\vbox to 8.5pt{}\voidright $}}}
\def\tenBig#1{{\hbox{$\left#1\vbox to 9.5pt{}\voidright $}}}
\def\tenbigg#1{{\hbox{$\left#1\vbox to 12.5pt{}\voidright $}}}
\def\tenBigg#1{{\hbox{$\left#1\vbox to 16pt{}\voidright $}}}
\def\elevenbig#1{{\hbox{$\left#1\vbox to 9pt{}\voidright $}}}
\def\elevenBig#1{{\hbox{$\left#1\vbox to 10.5pt{}\voidright $}}}
\def\elevenbigg#1{{\hbox{$\left#1\vbox to 14pt{}\voidright $}}}
\def\elevenBigg#1{{\hbox{$\left#1\vbox to 17.5pt{}\voidright $}}}
\def\twelvebig#1{{\hbox{$\left#1\vbox to 10pt{}\voidright $}}}
\def\twelveBig#1{{\hbox{$\left#1\vbox to 11pt{}\voidright $}}}
\def\twelvebigg#1{{\hbox{$\left#1\vbox to 15pt{}\voidright $}}}
\def\twelveBigg#1{{\hbox{$\left#1\vbox to 19pt{}\voidright $}}}
\def\fifteenbig#1{{\hbox{$\left#1\vbox to 12pt{}\voidright $}}}
\def\fifteenBig#1{{\hbox{$\left#1\vbox to 13.5pt{}\voidright $}}}
\def\fifteenbigg#1{{\hbox{$\left#1\vbox to 18pt{}\voidright $}}}
\def\fifteenBigg#1{{\hbox{$\left#1\vbox to 23pt{}\voidright $}}}
\def\voidright{\right.\nulldelimiterspace=0pt \mathsurround=0pt }
\def\fifteenpoint{
  \textfont0=\fifteenrm \scriptfont0=\twelverm \scriptscriptfont0=\tenrm
  \def\rm{\fam0 \fifteenrm}%
  \textfont1=\fifteeni \scriptfont1=\twelvei \scriptscriptfont1=\teni
  \textfont2=\fifteensy \scriptfont2=\twelvesy \scriptscriptfont2=\tensy
  \textfont3=\fifteenex \scriptfont3=\fifteenex \scriptscriptfont3=\fifteenex
  \def\it{\fam\itfam\fifteenit}\textfont\itfam=\fifteenit
  \def\sl{\fam\slfam\fifteensl}\textfont\slfam=\fifteensl
  \def\bf{\fam\bffam\fifteenbf}\textfont\bffam=\fifteenbf 
    \scriptfont\bffam=\twelvebf\scriptscriptfont\bffam=\tenbf
  \def\tt{\fam\ttfam\fifteentt}\textfont\ttfam=\fifteentt
  \def\sc{\fam\scfam\fifteensc}\textfont\scfam=\fifteensc
  \def\am{\fam\amfam\fifteenam}\textfont\amfam=\fifteenam
    \scriptfont\amfam=\twelveam\scriptscriptfont\amfam=\tenam
  \def\bm{\fam\bmfam\fifteenbm}\textfont\bmfam=\fifteenbm
    \scriptfont\bmfam=\twelvebm\scriptscriptfont\bmfam=\tenbm
  \baselineskip=21pt \rm
  \let\big=\fifteenbig\let\Big=\fifteenBig\let\bigg=\fifteenbigg
  \let\Bigg=\fifteenBigg}
\def\twelvepoint{
  \textfont0=\twelverm \scriptfont0=\ninerm \scriptscriptfont0=\sevenrm
  \def\rm{\fam0 \twelverm}%
  \textfont1=\twelvei \scriptfont1=\ninei \scriptscriptfont1=\seveni
  \textfont2=\twelvesy \scriptfont2=\ninesy \scriptscriptfont2=\sevensy
  \textfont3=\twelveex \scriptfont3=\twelveex \scriptscriptfont3=\twelveex
  \def\it{\fam\itfam\twelveit}\textfont\itfam=\twelveit
  \def\sl{\fam\slfam\twelvesl}\textfont\slfam=\twelvesl
  \def\bf{\fam\bffam\twelvebf}\textfont\bffam=\twelvebf 
    \scriptfont\bffam=\ninebf\scriptscriptfont\bffam=\sevenbf
  \def\tt{\fam\ttfam\twelvett}\textfont\ttfam=\twelvett
  \def\sc{\fam\scfam\twelvesc}\textfont\scfam=\twelvesc
  \def\am{\fam\amfam\twelveam}\textfont\amfam=\twelveam
    \scriptfont\amfam=\nineam\scriptscriptfont\amfam=\sevenam
  \def\bm{\fam\bmfam\twelvebm}\textfont\bmfam=\twelvebm
    \scriptfont\bmfam=\ninebm\scriptscriptfont\bmfam=\sevenbm
  \baselineskip=17.8pt \rm 
  \def\looselineskip{\baselineskip=18.5pt plus 1.8pt}%
  \def\tightlineskip{\baselineskip=16.5pt}%
  \def\verytightlineskip{\baselineskip=15pt}%
  \let\big=\twelvebig\let\Big=\twelveBig\let\bigg=\twelvebigg
  \let\Bigg=\twelveBigg  }
\def\elevenpoint{
  \textfont0=\elevenrm \scriptfont0=\ninerm \scriptscriptfont0=\sixrm
  \def\rm{\fam0 \elevenrm}%
  \textfont1=\eleveni \scriptfont1=\ninei \scriptscriptfont1=\sixi
  \textfont2=\elevensy \scriptfont2=\ninesy \scriptfont2=\sixsy 
  \textfont3=\elevenex \scriptfont3=\elevenex \scriptfont3=\elevenex
  \def\it{\fam\itfam\elevenit}\textfont\itfam=\elevenit
  \def\sl{\fam\slfam\elevensl}\textfont\slfam=\elevensl
  \def\bf{\fam\bffam\elevenbf}\textfont\bffam=\elevenbf
  \def\tt{\fam\ttfam\eleventt}\textfont\ttfam=\eleventt
  \def\sc{\fam\scfam\elevensc}\textfont\scfam=\elevensc
  \def\am{\fam\amfam\elevenam}\textfont\amfam=\elevenam
    \scriptfont\amfam=\nineam\scriptscriptfont\amfam=\sixam
  \def\bm{\fam\bmfam\elevenbm}\textfont\bmfam=\elevenbm
    \scriptfont\bmfam=\ninebm\scriptscriptfont\bmfam=\sixbm
  \baselineskip=15.1pt \rm
  \def\looselineskip{\baselineskip=16pt plus 1.5pt}%
  \def\tightlineskip{\baselineskip=14pt}%
  \def\verytightlineskip{\baselineskip=13pt}%
  \let\big=\elevenbig\let\Big=\elevenBig\let\bigg=\elevenbigg
  \let\Bigg=\elevenBigg  }
\def\tenpoint{
  \textfont0=\tenrm \scriptfont0=\eightrm \scriptscriptfont0=\fiverm
  \def\rm{\fam0 \tenrm}%
  \textfont1=\teni \scriptfont1=\eighti \scriptscriptfont1=\fivei
  \textfont2=\tensy \scriptfont2=\eightsy \scriptfont2=\fivesy 
  \textfont3=\tenex \scriptfont3=\tenex \scriptfont3=\tenex
  \def\it{\fam\itfam\tenit}\textfont\itfam=\tenit
  \def\sl{\fam\slfam\tensl}\textfont\slfam=\tensl
  \def\bf{\fam\bffam\tenbf}\textfont\bffam=\tenbf
  \def\tt{\fam\ttfam\tentt}\textfont\ttfam=\tentt
  \def\sc{\fam\scfam\tensc}\textfont\scfam=\tensc
  \def\am{\fam\amfam\tenam}\textfont\amfam=\tenam
    \scriptfont\amfam=\eightam \scriptscriptfont\amfam=\fiveam
  \def\bm{\fam\bmfam\tenbm}\textfont\bmfam=\tenbm
    \scriptfont\bmfam=\eightbm \scriptscriptfont\bmfam=\fivebm
  \baselineskip=14pt \rm
  \def\looselineskip{\baselineskip=14.8pt plus1.5pt}
  \def\tightlineskip{\baselineskip=13.6pt}%
  \def\verytightlineskip{\baselineskip=13pt}%
  \let\big=\tenbig\let\Big=\tenBig\let\bigg=\tenbigg\let\Bigg=\tenBigg  }
\def\ninepoint{
  \textfont0=\ninerm \scriptfont0=\sevenrm \scriptscriptfont0=\fourrm
  \def\rm{\fam0 \ninerm}%
  \textfont1=\ninei \scriptfont1=\seveni \scriptscriptfont1=\fouri
  \textfont2=\ninesy \scriptfont2=\sevensy \scriptfont2=\foursy 
  \textfont3=\nineex \scriptfont3=\nineex \scriptfont3=\nineex
  \def\it{\fam\itfam\nineit}\textfont\itfam=\nineit
  \def\sl{\fam\slfam\ninesl}\textfont\slfam=\ninesl
  \def\bf{\fam\bffam\ninebf}\textfont\bffam=\ninebf
  \def\tt{\fam\ttfam\ninett}\textfont\ttfam=\ninett
  \def\sc{\fam\scfam\ninesc}\textfont\scfam=\ninesc
  \def\am{\fam\amfam\nineam}\textfont\amfam=\nineam
    \scriptfont\amfam=\nineam\scriptscriptfont\amfam=\fouram
  \def\bm{\fam\bmfam\ninebm}\textfont\bmfam=\ninebm
    \scriptfont\bmfam=\ninebm\scriptscriptfont\bmfam=\fourbm
  \baselineskip=12.6pt \rm
  \let\big=\ninebig\let\Big=\nineBig\let\bigg=\ninebigg
  \let\Bigg=\nineBigg  }
\def\eightpoint{
  \textfont0=\eightrm \scriptfont0=\fiverm \scriptscriptfont0=\fourrm
  \def\rm{\fam0 \eightrm}%
  \textfont1=\eighti \scriptfont1=\fivei \scriptscriptfont1=\fouri
  \textfont2=\eightsy \scriptfont2=\fivesy \scriptfont2=\foursy 
  \textfont3=\eightex \scriptfont3=\eightex \scriptfont3=\eightex
  \def\it{\fam\itfam\eightit}\textfont\itfam=\eightit
  \def\sl{\fam\slfam\eightsl}\textfont\slfam=\eightsl
  \def\bf{\fam\bffam\eightbf}\textfont\bffam=\eightbf
  \def\tt{\fam\ttfam\eighttt}\textfont\ttfam=\eighttt
  \def\sc{\fam\scfam\eightsc}\textfont\scfam=\eightsc
  \def\am{\fam\amfam\eightam}\textfont\amfam=\eightam
    \scriptfont\amfam=\eightam\scriptscriptfont\amfam=\fouram
  \def\bm{\fam\bmfam\eightbm}\textfont\bmfam=\eightbm
    \scriptfont\bmfam=\eightbm\scriptscriptfont\bmfam=\fourbm
  \baselineskip=11.2pt \rm
  \let\big=\eightbig\let\Big=\eightBig\let\bigg=\eightbigg
  \let\Bigg=\eightBigg  }

\twelvepoint
\nopagenumbers
\hsize=6in\vsize=8.8in

\parskip=1pt plus 1pt

\newif\ifSpecialhead\Specialheadfalse
\newbox\specialheadbox

\def\specialhead #1\par{\Specialheadtrue\setbox\specialheadbox=\hbox{#1}}
\headline={{\ifSpecialhead\box\specialheadbox\global\Specialheadfalse\else
     \ifnum\pageno<0{\hfill\quad{\twelvebf\folio}}%
     \else\ifnum\pageno<2\hfill
     \else\hfill\twelvepoint\sc\firstmark\quad{\twelvebf\folio}\fi\fi\fi}}

\def\title#1\par{\bigskip{\def\cr{\par\center}\center\fifteenbf #1\par}\medskip}
\def\subtitle#1\par{\centerline{\fifteenrm #1}\medskip}
\def\author#1\par{\medskip{\def\cr{\par\center\twelvesc}\fifteensc\center#1\par}}
\def\center#1\par{\hfil #1\hfil\par}
\def\abstract.#1\par{\message{Abstract.}%
                    \medskip{\narrower\narrower\tenpoint\tightlineskip
                        \noindent{\bf Abstract.}#1\par}\medskip\noindent}
\def\bigabstract.#1\par{\message{Abstract.}%
                         \medskip{\narrower\narrower\tightlineskip
                         \noindent{\bf Abstract. }#1\par}\medskip\noindent}
\def\acknowledgement#1\par{\footnote{}{#1}}
\def\sectionskip{\Goodbreak\vskip 25pt plus 15pt minus 5pt}
\def\secnumber{\ifquiet
               \else\ifNoSections
                    \else\sectionsymbol\the\secno\quad\fi\fi}
\def\section#1\par{ \NoSectionsfalse\par\sectionskip\proofdepth=0\claimno=0
 \ifquiet\else\advance\secno by1\fi\toks0={#1}
 \immediate\write16{\ifquiet\else Section \the\secno\space\fi
                    \the\toks0}%
 \mark{\secnumber #1}%
 {\fifteenpoint\bf\noindent\secnumber #1}\nobreak\bigskip\quietoff
 \nobreak\noindent}

\def\QUIET{\QUIETtrue\quiettrue}

\def\quietoff{\ifQUIET\else\quietfalse\fi}
\newif\ifquiet
\newif\ifQUIET
\newif\ifNoSections
\newcount\claimtype
\newcount\secno
\newcount\claimno
\newcount\subclaimno
\newcount\subsubclaimno
\newcount\subsubsubclaimno
\newcount\proofdepth
\def\subclaimnumber{\ifquiet\else\ifcase\subclaimno\or A\or B\or C\or D\or E\or
     F\or G\or H\or I\or J\or K\or L\or M\or N\or O\or P\fi\fi}
\def\subsubclaimnumber{\ifquiet\else\ifcase\subsubclaimno\or i\or ii\or iii\or 
   iv\or v\or vi\or vii\or viii\or ix\or x\or xi\or xii\or xiii\or xiv\fi\fi}
\def\subsubsubclaimnumber{\ifquiet\else\ifcase\subsubsubclaimno\or a\or b\or 
   c\or d\or e\or f\or g\or \or h\or i\or j\or k\or l\or m\or n\or o\fi\fi}
\def\claimtag{\ifquiet\else
  \ifNoSections
    \ifcase\proofdepth\the\claimno%
    \or\the\claimno.\subclaimnumber
    \or\the\claimno.\subclaimnumber.\subsubclaimnumber
    \or\the\claimno.\subclaimnumber.\subsubclaimnumber
                                                .\subsubsubclaimnumber\fi
  \else
    \ifcase\proofdepth\the\secno.\the\claimno
    \or\the\secno.\the\claimno.\subclaimnumber
    \or\the\secno.\the\claimno.\subclaimnumber.\subsubclaimnumber
    \or\the\secno.\the\claimno.\subclaimnumber.\subsubclaimnumber
                                                .\subsubsubclaimnumber\fi\fi\fi}
\secno=0\claimno=0\proofdepth=0\subclaimno=0\subsubclaimno=0\subsubsubclaimno=0
\NoSectionstrue
\newbox\qedbox
\def\claimname{\ifcase\claimtype Theorem\or Lemma\or Claim\or Corollary\or
               Question\or Definition\or Remark\or Conjecture\fi}
\def\preclaimskip{\removelastskip
    \ifcase\claimtype\goodbreak\vskip 8pt plus 4pt minus 2pt
                  \or\goodbreak\vskip 6pt plus 4pt minus 1pt
                  \or\goodbreak\vskip 5pt plus 4pt minus 1pt
                  \or\goodbreak\vskip 8pt plus 4pt minus 2pt
                  \or\vskip 7pt plus 4pt minus 2pt
                  \or\vskip 7pt plus 4pt minus 2pt
                  \or\vskip 7pt plus 4pt minus 2pt
                  \or\goodbreak\vskip 8pt plus 4pt minus 2pt\fi}
\def\postclaimskip{\ifcase\claimtype         \vskip 4pt plus 2pt minus 2pt
                                          \or\vskip 3pt plus 2pt minus 2pt
                                          \or\vskip 2pt plus 2pt minus 1pt
                                          \or\vskip 4pt plus 2pt minus 2pt
                                          \or\vskip 1pt plus 2pt 
                                          \or\vskip 4pt plus 4pt 
                                          \or\vskip 3pt plus 2pt
                                          \or\vskip 4pt plus 2pt minus 2pt\fi}
\def\claimfont{\ifcase\claimtype
                  \sl\or\sl\or\sl\or\sl\or\sl\or\rm\or\rm\or\sl\fi}
\def\advancetag{\ifcase\proofdepth\advance\claimno by1
                               \or\advance\subclaimno by1
                               \or\advance\subsubclaimno by1
                               \or\advance\subsubsubclaimno by1\fi}
\def\sayclaim#1.#2 #3\par{\ifquiet\else\advancetag\fi
    \preclaimskip\setbox1=\hbox{#1}\setbox2=\hbox{#2}%
    \toks0={#1 }
    \immediate\write16{\ifdim\wd1>0pt\the\toks0
                       \else\claimname\space\fi \claimtag.}%
    \vbox{\noindent
    {\bf\ifdim\wd1=0pt \claimname\else #1\fi\ifquiet.\else\ \claimtag{\ifNoSections.\fi}\fi}%
    \enspace{\ifdim\wd2>0pt\sc #2\enspace\fi}%
    {\claimfont #3\par}}\postclaimskip\quietoff}
\def\theorem{\claimtype=0\sayclaim}
\def\lemma{\claimtype=1\sayclaim}

\def\question{\claimtype=4\sayclaim}

\def\remark{\claimtype=6\sayclaim}

\def\point#1. #2\par{\item{\rm #1.}#2\par}
\def\points#1\cr\par{\medskip\vbox{\let\cr=\point\point#1\par}\par}

\def\prooffont{}
\def\proofsize{}
\def\proofindent{}
\def\proofskip{\badbreak\ifcase\claimtype    \vskip 3pt plus 2pt minus 2pt
                                          \or\vskip 2pt plus 2pt minus 2pt
                                          \or\vskip 1pt plus 2pt minus 1pt
                                          \or\vskip 3pt plus 2pt minus 2pt
                                          \or\vskip 1pt plus 2pt 
                                          \or\vskip 2pt plus 4pt 
                                          \or\vskip 1pt plus 2pt
                                          \or\vskip 3pt plus 2pt minus 2pt\fi}
\def\greatbreak{\vskip0pt plus1.5in\penalty-1000\vskip0pt plus-2in}

\def\Goodbreak{\vskip0pt plus.5in\penalty-1000\vskip0pt plus-.5in}
\def\goodbreak{\penalty-500}
\def\badbreak{\penalty500}
\def\Badbreak{\penalty1000}
\def\proof{\message{proof}\removelastskip\Badbreak\proofskip\begingroup
  \advance\proofdepth by1
  \setbox\qedbox=\hbox{\halmos\raise2pt\hbox{\fiverm\claimname}}%
  \prooffont\proofsize\proofindent\noindent{\bf Proof: }}
\def\proofof#1:{\message{proof}\removelastskip\Badbreak\proofskip\begingroup
  \advance\proofdepth by1
  \setbox\qedbox=\hbox{\halmos\raise2pt\hbox{\fiverm#1}}%
  \prooffont\proofsize\proofindent\noindent{\bf Proof of #1: }}
\def\cite[#1]{[{\tenrm{#1}}]\message{[#1]}}
\edef\ref#1{\expandafter\global\expandafter\edef#1{\noexpand\claimtag}}
\newwrite\notes
\openout\notes=\jobname.notes
\long\def\unexpandedwrite#1#2{\def\finwrite{\write#1}%
   {\aftergroup\finwrite\aftergroup{\sanitize#2\endsanity}}}
\def\sanitize{\futurelet\next\sanswitch}
\let\stoken=\space
\def\sanswitch{\ifx\next\endsanity
  \else\ifcat\noexpand\next\stoken\aftergroup\space\let\next=\eat
   \else\ifcat\noexpand\next\bgroup\aftergroup{\let\next=\eat
    \else\ifcat\noexpand\next\egroup\aftergroup}\let\next=\eat
     \else\let\next=\copytoken\fi\fi\fi\fi \next}
\def\eat{\afterassignment\sanitize \let\next= }
\long\def\copytoken#1{\ifcat\noexpand#1\relax\aftergroup\noexpand
  \else\ifcat\noexpand#1\noexpand~\aftergroup\noexpand\fi\fi
  \aftergroup#1\sanitize}
\def\endsanity\endsanity{}

\def\note#1#2{\hbox to2in{\strut#1\quad\dotfill\quad#2}}
\def\boxit#1{\setbox4=\hbox{\kern1pt#1\kern1pt}
  \hbox{\vrule\vbox{\hrule\kern1pt\box4\kern1pt\hrule}\vrule}}
\def\halmos{\hbox{\am\char'3}} 
\def\qed#1\par{\message{.                                }\setbox1=\hbox{#1}%
  \ifdim\wd1>0pt\setbox\qedbox=\hbox{\halmos\raise2pt\hbox{\fiverm #1}}\fi
  \kern5pt\lower 2pt\hbox{\box\qedbox}\proofskip\goodbreak\endgroup}

\def\sectionsymbol{\S}
\def\k{\kappa}
\def\g{\gamma}
\def\a{\alpha}
\def\b{\beta}
\def\d{\delta}

\def\l{\lambda}

\def\I1{\mathop{\hbox{\sc i}_1}}
\def\w{\omega}

\def\Z{{\mathchoice{\hbox{\bm Z}}{\hbox{\bm Z}}
         {\hbox{\tenbm Z}}{\hbox{\sevenbm Z}}}}
\def\N{{\mathchoice{\hbox{\bm N}}{\hbox{\bm N}}
         {\hbox{\tenbm N}}{\hbox{\sevenbm N}}}}

\def\card#1{\left|#1\right|}

\def\dirlim{\mathop{\rm dir\,lim}\nolimits}

\def\aut{\mathop{\hbox{\rm Aut}}\nolimits}
\def\inn{\mathop{\hbox{\rm Inn}}\nolimits}

\def\unifto{\buildrel\lower 7pt\hbox{$\to$}\over\to}

\def\iso{\cong}

\def\cof{\mathop{\rm cof}\nolimits}

\def\ZFC{\hbox{\sc zfc}}
\def\GCH{\hbox{\sc gch}}

\def\plus{^{\scriptscriptstyle +}}

\def\in{\mathrel{\mathchoice{\raise 
1pt\hbox{$\scriptstyle\cal\char'62$}}
         {\raise 1pt\hbox{$\scriptstyle\cal\char'62$}}
         {\raise .5pt\hbox{$\scriptscriptstyle\cal\char'62$}}
         {\hbox{$\scriptscriptstyle\cal\char'62$}}}\penalty700{}}
\def\ni{\mathrel{\mathchoice{\raise 1pt\hbox{$\scriptstyle\cal\char'63$}}
                   {\raise 1pt\hbox{$\scriptstyle\cal\char'63$}}
                   {\raise .5pt\hbox{$\scriptscriptstyle\cal\char'63$}}
                   {\hbox{$\scriptscriptstyle\cal\char'63$}}}\penalty700}
\def\of{\mathrel{\mathchoice{\raise 1pt\hbox{$\scriptstyle\subseteq$}}
                   {\raise 1pt\hbox{$\scriptstyle\subseteq$}}
                   {\raise .5pt\hbox{$\scriptscriptstyle\subseteq$}}
                   {\hbox{$\scriptscriptstyle\subseteq$}}}}
\def\fo{\mathrel{\mathchoice{\raise 1pt\hbox{$\scriptstyle\supseteq$}}
                   {\raise 1pt\hbox{$\scriptstyle\supseteq$}}
                   {\raise .5pt\hbox{$\scriptscriptstyle\supseteq$}}
                   {\hbox{$\scriptscriptstyle\supseteq$}}}}
\def\notin{\mathrel{\mathchoice
  {\raise 1pt\hbox{\rlap{$\scriptstyle\;|$}$\scriptstyle\cal\char'62$}}
  {\raise 1pt\hbox{\rlap{$\scriptstyle\kern2pt 
          |$}$\scriptstyle\cal\char'62$}}
  {\raise .5pt\hbox{\rlap{$\scriptscriptstyle\, |$}$\scriptscriptstyle
      \cal\char'62$}}
  {\hbox{\rlap{$\scriptscriptstyle\, |$}$\scriptscriptstyle
     \cal\char'62$}}}%
  \penalty700}

\def\and{\mathrel{\kern1pt\&\kern1pt}}

\def\union{\cup}

\def\compose{\circ}

\def\cross{\times}

\def\tlt{\triangleleft}

\def\[#1]{\left[\vphantom{\bigm|}#1\right]}
\def\<#1>{\langle\,#1\,\rangle}

\def\restrict{\mathbin{\mathchoice{\hbox{\am\char'26}}{\hbox{\am\char'26}}{\hbox{\eightam\char'26}}{\hbox{\sixam\char'26}}}}

\def\st{\mid}
\def\seq<#1>{{\def\st{\mid\penalty650}\left<\,#1\,\right>}}

\def\set#1{\{\,#1\,\}}

\def\lttheta{{\raise 1pt\hbox{$\scriptstyle<$}\theta}}

\def\I1{\mathop{\hbox{\sc i}_1}}

\def\Diamond{\diamondsuit}

\QUIET

\center [Talk for the Mathematical Society of Japan in Osaka on October 2, 1998]

\title How tall is the automorphism tower of a group?

\author Joel David Hamkins\cr
            Kobe University and\cr
	The City University of New York\cr
	{\ninett http://www.library.csi.cuny.edu/users/hamkins}\cr

\noindent The automorphism tower of a group is obtained by computing
its automorphism group, the automorphism group of {\it that} group, and
so on, iterating transfinitely. Each group maps canonically into the
next using inner automorphisms, and so at limit stages one can take a
direct limit and continue the iteration.  $$G\to\aut(G)\to
\aut(\aut(G))\to\cdots\to G_\omega\to G_{\omega+1}\to \cdots\to
G_\a\to\cdots$$ The tower is said to terminate if a fixed point is
reached, that is, if a group is reached which is isomorphic to its
automorphism group by the natural map. This occurs if a {\it complete}
group is reached, one which is centerless and has only inner
automorphisms.

The natural map $\pi:G\to\aut(G)$ is the one that takes any element $g\in G$ to the inner automorphism $i_g$, defined by simple conjugation $i_g(h)=ghg^{-1}$. Thus, the kernel of $\pi$ is precisely the center of $G$, the set of elements which commute with everything in $G$, and the range of $\pi$ is precisely the set of inner automorphisms of $G$. By composing the natural maps at every step, one obtains a commuting system of homomorphisms $\pi_{\a,\b}:G_\a\to G_\b$ for $\a<\b$, and these are the maps which are used to compute the direct limit at limit stages. 

Much of the historical analysis of the automorphism tower has focused on the special case of centerless groups, for when the initial group is centerless, matters simplify considerably\footnote{${}^*$}{In Hulse [1970], Rae and
Roseblade [1970], and Thomas [1985], the tower is only defined in this
special case; but the definition I give here works perfectly well
whether or not the group $G$ is centerless. Of course, when there is a
center, one has homomorphisms rather than embeddings.}. An easy computation shows that $\theta\compose i_g\compose\theta=i_{\theta(g)}$ for any automorphism $\theta$, and from this we conclude that $\inn(G)\tlt\aut(G)$ and, for centerless $G$, that $C_{\aut(G)}(\inn(G))=1$. In particular, if $G$ is centerless then so also is $\aut(G)$, and more generally, by tranfinite induction every group in the automorphism tower of a centerless group is centerless. In this case, consequently, all the natural maps $\pi_{\a,\b}$ are injective, and so by identifying every group with its image under the canonical map, we may view the tower as building upwards to larger and larger groups; the question is whether this building process ever stops.
$$G\of G_1\of \cdots G_\w\of\cdots\of G_\a\of\cdots$$ In this centerless case, the outer automorphisms of every group become inner automorphisms in the next group. One wants to know, then, whether this process eventually closes off. 

The classical result is the following theorem of Wielandt.

\theorem Classical Theorem.{(Wielandt, 1939)} The automorphism tower of any centerless finite group terminates in finitely many steps. 

Wielandt's theorem, pointed to with admiration at the conclusion of Scott's [1964] book {\it Group Theory}, was the inspiration for a line of gradual generalizations by various mathematicians over the succeeding decades. Scott closes his book with questions concerning the automorphism tower, specifically mentioning the possibility of transfinite iterations and towers of non-centerless groups.

\question.{(Scott, 1964)} Is there a group whose automorphism tower never terminates?

By the 1970's, several authors had made progress: 

\theorem.{(Rae and Roseblade, 1970)} The automorphism tower of any centerless \v Cernikov group terminates in finitely many steps. 

\theorem.{(Hulse 1970)} The automorphism tower of any centerless polycyclic group terminates in a countable ordinal number of steps. 

These results culminated in Simon Thomas' elegant solution to the automorphism tower problem in the case of centerless groups. An application of Fodor's lemma lies at the heart of Thomas' proof. 

\theorem.{(Thomas, 1985)} The automorphism tower of any centerless group eventually terminates. Indeed, the automorphism tower of a centerless group $G$ terminates in fewer than $(2^{\card{G}})\plus$ many steps.  

In the general case, however, the question remained open whether every group has a terminating automorphism tower.  This is settled by the following theorem.

\theorem Main Theorem.{(Hamkins [1998])} Every group has a terminating automorphism tower.

\proof Suppose $G$ is a group. The following transfinite 
recursion defines the automorphism tower of $G$:
$$\eqalign{G_0=&\,G\cr
           G_{\a+1}=&\aut(G_\a), \quad
           \hbox{where $\pi_{\a,\a+1}:G_\a\to G_{\a+1}$ is the natural map,}\cr
           G_\l=&\dirlim_{\a<\l}G_\a, \quad
            \hbox{if $\l$ is a limit ordinal.}\cr}$$
When $\a<\b$ one obtains the map $\pi_{\a,\b}:G_\a\to G_\b$ by
composing the canonical maps at each step, and these are the maps used
to compute the direct limit at limit stages. Thus, when $\l$ is a limit
ordinal, every element of $G_\l$ is of the form $\pi_{\a,\l}(g)$ for
some $\a<\l$ and some $g\in G_\a$.

Since Simon Thomas [1985] has proved that every centerless group has a
terminating automorphism tower, it suffices to show that there is an
ordinal $\g$ such that $G_\g$ has a trivial center. For each ordinal
$\a$, let $H_\a=\set{g\in G_\a\st \exists\b\,\pi_{\a,\b}(g)=1}$.  For
every $g\in H_\a$ there is some least $\b_g>\a$ such that
$\pi_{\a,\b_g}(g)=1$.  Let $f(\a)=\sup_{g\in H_\a}\b_g$. It is easy to
check that if $\a<\b$ then $\a<f(\a)\leq f(\b)$. Iterating the
function, define $\g_0=0$ and $\g_{n+1}=f(\g_n)$. This produces a
strictly increasing $\w$-sequence of ordinals whose supremum
$\g=\sup\set{\g_n\st n\in\w}$ is a limit ordinal which is closed under
$f$. That is, $f(\a)<\g$ for every $\a<\g$.  I claim that $G_\g$ has a
trivial center. To see this, suppose $g$ is in the center of $G_\g$.
Thus, $\pi_{\g,\g+1}(g)=1$. Moreover, since $\g$ is a limit ordinal,
there is $\a<\g$ and $h\in G_\a$ such that $g=\pi_{\a,\g}(h)$.
Combining these facts, observe that
$$\pi_{\a,\g+1}(h)=\pi_{\g,\g+1}(\pi_{\a,\g}(h))=\pi_{\g,\g+1}(g)=1.$$
Consequently, $\pi_{\a,f(\a)}(h)=1$.  Since $f(\a)<\g$, it follows that
$$g=\pi_{\a,\g}(h)=\pi_{f(\a),\g}(\pi_{\a,f(\a)}(h))=\pi_{f(\a),\g}(1)=1,$$
as desired.\qed

So now we know that the automorphism tower of any group terminates. But how long does it take? In the centerless case, Simon Thomas provided an attractive bound on the height of the automorphism tower, namely, the tower of a centerless group $G$ terminates before $(2^{\card{G}})\plus$. It is therefore natural to ask the question for groups in general:

\question. How tall is the automorphism tower of a group $G$?

Unfortunately, the proof of the Main Theorem above does not reveal exactly how long the automorphism tower takes
to stabilize, since it is not clear how large $f(\a)$ can be. Nevertheless, something more can be said. Certainly the automorphism
tower of $G$ terminates well before the next inaccessible cardinal
above $\card{G}$.  More generally, if $\l>\w$ is regular and
$\card{G_\a}<\l$ whenever $\a<\l$, then I claim the centerless groups
will appear before $\l$. To see this, let $H_\a=\set{g\in G_\a\st
\exists\b_g{<}\l\,\, \pi_{\a,\b_g}(g)=1}$ and define $f:\l\to\l$ by
$f(\a)=\sup_{g\in H_\a}\b_g$; if $\g<\l$ is closed under $f$, then it
follows as in the main theorem that $G_\g$ has no center, as desired.
In this case the tower therefore terminates in fewer than
$(2^{\card{G_\g}})\plus$ many additional steps. If it happens that
$\card{G_\g}\plus<\l$, one can adapt Thomas' \cite[1996] argument using
Fodor's lemma to prove that the tower terminates actually in fewer than
$\l$ many steps. The point is that one can find a bound on the height
of the tower by bounding the rate of growth of the groups in the
tower.

Thomas \cite[1985] provides his explicit bound in the case of
centerless $G$ in precisely this way. He proves that if $G$ is
centerless and $\l=(2^{\card{G}})\plus$, then $\card{G_\a}<\l$ for all
$\a<\l$. The analogous result, unfortunately, does not hold for groups
with nontrivial centers. This is illustrated by the following example,
provided by the anonymous referee of  my paper [1998]:

\theorem Example. There exists a countable group $G$ such that 
$\card{\aut G}=2^\w$ and\break $\card{\aut(\aut G)}=2^{2^\w}$.

\proof For each prime $p$, let $\Z[1/p]=\set{m/{p^n}\st m\in\Z,n\in\N}$
be the additive group of $p$-adic rationals and let $G=\oplus_p\Z[1/p]$
be the direct sum of these groups. An element $g\in G$ is divisible by $p^n$
for all $n\in\N$ iff $g\in\Z[1/p]$. Hence, if $\pi\in\aut G$, then
$\pi[\Z[1/p]]=\Z[1/p]$ for each prime $p$. It is easy to see that any
automorphism of $\Z[1/p]$ is simply multiplication by an element $u\in
U_p=\set{\pm p^n\st n\in\Z}$, the group of multiplicative units of the
ring of $p$-adic rationals. Thus, $\aut G\iso\Pi_p U_p$; and so
$\card{\aut G}=2^\w$.

Next note that $U_p\iso\Z\cross C_2$ for each prime $p$. Thus $\aut
G\iso P\cross V$, where $P$ is the direct product of countably many
copies of $\Z$ and $V$ is the direct product of countably many copies
of $C_2$. Since each nonzero element of $V$ has order $2$, it follows
that $V$ is isomorphic to a direct sum of $\card{V}$ copies of $C_2$.
Thus we can identify $V$ with a vector space of dimension $2^\w$ over
the field of two elements. Hence, $$\aut(\aut G)\iso \aut P\cross\aut
V=\aut P\cross GL(V),$$ where $GL(V)$ is the general linear group on
the vector space $V$. Since $\card{GL(V)}=2^{2^\w}$, it follows that
$\card{\aut(\aut G)}=2^{2^\w}$.\qed{ }

Enriqueta Rodr\'\i guez-Carrington has observed that the natural
modification of my argument shows that the derivation tower of every
Lie algebra eventually leads to a centerless Lie algebra. Since Simon
Thomas [1985] proved that the derivation tower of every centerless Lie
algebra must eventually terminate, it follows that the derivation tower
of any Lie algebra must eventually terminate.

One might hope, since every step of the automorphism tower kills the
center of the previous group, that $G_\omega$ is always centerless; but
this is not so. The dihedral group with eight elements has a center of
size two, but is isomorphic to its own automorphism group (there is an
outer automorphism which swaps $a$ and $b$ in the presentation
$\<a,b\st a^2=1,b^2=1,(ab)^4=1>$ ).  The group at stage $\omega$ is
just the two element group, which still has a center, and so this tower
survives until $\omega+1$. Simon Thomas has constructed examples
showing that for every natural number $n$ there are finite groups whose
tower has height $\omega+n$, but these also become centerless at stage
$\w+1$. Perhaps our attention should focus, for an arbitrary group $G$,
on the least ordinal stage $\g$ such that $G_\g$ is centerless. The
main point, then, is to find an explicit bound on how large $\g$ can be
in comparison with $\card{G}$.

Let me now ask an innocent question:

\question. Can you predict the height of the automorphism tower of a group $G$ by looking at $G$?

Of course, you may argue philosophically, the answer is Yes, becasue the automorphism tower of a group $G$ is completely determined by $G$; one simply iterates the automorphism group operation until the termination point is obtained, and that is where the tower terminates. But nevertheless, I counter philosophically that the answer to the question is No! How can this be? 

The reason for my negative answer is that the automorphism tower of a group has a set-theoretic essence; building the automorphism tower by iteratively computing automorphism groups is rather like building the Levy hierarchy $V_\a$ by iteratively computing power sets. The fact is that the automorphism tower of a group can depend on the model of set theory in which you compute it, with the very same group leading to wildly different automorphism towers in different set theoretic universes. And so, in order to predict the height of the tower of $G$, you can't just look at $G$; you must also look at $\aut(G)$ and $\aut(\aut(G))$ and so on; and these depend on the set-theoretic background. 

\theorem.{(Hamkins and Thomas 1997)} It is relatively consistent that for every $\l$ and for every $\a<\l$ there is a group $G$ whose tower has height $\a$, but for any non-zero $\b<\l$ there is a forcing extension in which the height of the tower of $G$ has height $\b$. 

I would like for most of the remaining time to give a stratospheric view of the proof of this theorem. For the details, which are abundant, I refer you to our paper [1997].

\center {\tenbf The set-theoretic kernel}

\noindent Perhaps the key set-theoretic idea is the realization that forcing can make non-isomorphic structures isomorphic. We begin with the following problem.

\remark Warm-up Problem. Construct rigid non-isomorphic objects $S$ and $T$ which can be made isomorphic by forcing while remaining rigid. 

The solution is to use generic Souslin trees. If you add mutually generic normal Souslin trees, by forcing with normal $\a$-trees ($\a<\w_1$) ordered by end-extension, the resulting $\w_1$-trees $S$ and $T$ will be rigid and non-isomorphic. To see why this is so, suppose, for example, that $\dot\pi$ is the name of a purported isomorphism between $S$ and $T$. A simple bootstrap argument provides conditions $S_\a$ and $T_\a$ which decide $\dot\pi\restrict S_\a$, and by strategically extending some paths in $S_\a$ but not the image of the paths in $T_\a$, one obtains a stronger condition $\<S_{\a+1}, T_{\a+1}>$ which forces that $\dot\pi$ is not an isomorphism, a contradiction. A similar bootstrap argument shows that $S$ and $T$, individually, are rigid. Since by a back-and-forth argument any two normal $\a$-trees, for $\a<\w_1$, are isomorphic, there are abundant partial isomorphisms on the intial segments of $S$ and $T$, and by forcing with these partial isomorphisms, one adds by forcing an isomorphism between $S$ and $T$. The combined forcing can be viewed as forcing with normal $\a$-trees $S_\a$ and $T_\a$ and an isomorphism between them, and consequently the bootstrap argument can be modified to show that the isomorphism which is added between $S$ and $T$ is unique. That is, the trees $S$ and $T$ remain rigid even after they are forced to be isomorphic, as desired. 

More generally, by an Easton support forcing iteration, we obtain such objects simultaneously for every regular cardinal:

\theorem.{(Hamkins and Thomas 1997)} One can add, for every regular cardinal $\k$, a set $\set{T_\a\st\a<\k}$ of pair-wise non-isomorphic rigid trees such that for any equivalence relation $E$ on $\k\plus$, there is a forcing extension preserving rigidity in which the isomorphism relation on the trees is exactly $E$. 

By this theorem, therefore, we have a collection of non-isomorphic rigid objects which we can make isomorphic at will, while preserving their rigidity. Later, these objects will be the unit elements in elaborate graphs whose automorphism groups we want to control precisely by forcing. 

\greatbreak
\center {\tenbf The algebraic kernel}

\noindent We begin with the admission that automorphism towers are too difficult to handle, and so instead we work with the {\it normalizer tower} of a group. Given $G\leq H$ define 
$$\eqalign{N_0(G)&=G\cr
                  N_{\a+1}(G)&=N_H(N_\a(G))\cr
                  N_\l(G)&=\union_{\a<\l}N_\a(G), \hbox{if $\l$ is a limit}\cr}$$

\noindent The reason for doing so lies in the following amazing fact:

\lemma Fact. The automorphism tower of a centerless group $G$ is exactly the normalizer tower of $G$ computed in the terminal group $H=G_\g$. That is, $G_\a=N_\a(G)$. 

\lemma Fact. Conversely, if $H\leq\aut(K)$, where $K$ is a field, then by making a few modifications (adding a few bells and whistles) to the normalizer tower of $H$ in $\aut(K)$, one obtains an automorphism tower of the same height. 

Thus, to make automorphism towers of a specific height, we need only make normalizer towers of that height in the automorphism group of a field. This latter limitation is considerably loosened in light of the following theorem of Fried and Kollar:

\theorem.{(Fried and Kollar 1981)} By adding points, any graph can be made into a field with the same automorphism group. 

\noindent Thus, we can restrict our attention to subgroups of the automorphism groups of graphs. Since any tree can be represented as a graph, we make the connection with our earlier set-theoretic argument. There, we obtained a delicate skill to make trees isomorphic while preserving there rigidity. By combining these trees in elaborate combinations, we will construct graphs whose automorphism groups we can precisely modify by forcing. Let us see how this is done.

\def\tri{\triangle}
\def\square{\halmos}
\def\circle{\bigcirc}
\def\sp{\hskip 5pt}
\def\cboxit#1{\setbox4=\hbox{\kern2pt#1\kern2pt}
  \hbox{\vrule\vtop{\vbox{\hrule\kern1pt\box4}\kern1pt\hrule}\vrule}}
\def\boxtwotri{\cboxit{\sp$\tri\sp\tri$\sp}}
\def\boxfourtri{\cboxit{\sp\boxtwotri\sp\boxtwotri\sp}}
In order to construct the tall normalizer towers, we take as a unit some rigid graph, which I will denote by $\tri$, and build the following large graph and the corresponding subgroup of its automorphism group:
$$\tri\sp\tri\sp\boxtwotri\sp\boxfourtri\sp\cdots\cdots\hbox{up to $\a$}$$
The intended subgroup, a large wreath product, is indicated by the boxes. The idea is that while the full automorphism group can freely permute the $\tri$s, the subgroup we are interested in consists of those permutations which iteratively swap the two components of any of the boxes. Thus, every element of the subgroup fixes the first two $\tri$s, but there is, for example, a group element which swaps the third and the fourth $\tri$ and swaps the fifth and the sixth $\tri$ for the eighth and the seventh $\tri$, respectively. The subgroup is therefore simply a large wreath product. 

The normalizer tower of this group in the full automorophism group of the graph can be iteratively computed with ease. At the first step, one must add the permutations which swap the first two triangles, leading to:
$$\boxtwotri\sp\boxtwotri\sp\boxfourtri\sp\cdots\cdots\sp\hbox{up to $\a$}$$
At the next stage one adds the permutations which swap the first two boxes here, producing:
$$\boxfourtri\sp\boxfourtri\sp\cdots\cdots\sp\hbox{up to $\a$}$$
Inductively, one sees that the normalizer tower continues to grow for $\a$ many steps and then terminates. That is, the normalizer tower of the initial group inside the full automorphism group has height $\a$. By adding the bells and whistles I mentioned earlier, then, we have constructed a centerless group whose automorphism tower has height $\a$. 

Now we come to the exciting twist, which is to use the rigid pairwise non-isomorphic objects obtained from the set-theoretic part of the argument. Representing these various objects by $\tri$, $\square$, $\circle$ and $\Diamond$ we build a subgroup of the following graph:
$$\tri\sp\square\sp\cboxit{$\sp\circle\sp\circle\sp$}\sp\cboxit{\sp\cboxit{$\sp\Diamond\sp
\Diamond\sp$}\cboxit{$\sp\Diamond\sp\Diamond\sp$}\sp}\sp\cdots\cdots$$
The key point is now that by forcing the objects to be isomorphic $\tri\iso\square\iso\circle\iso\Diamond\iso\cdots$ up to $\b$, one thereby transforms the previous graph to look like the original picture 
$$\tri\sp\tri\sp\boxtwotri\sp\boxfourtri\sp\cdots\cdots$$
for $\b$ many steps. After the forcing, therefore, the normalizer tower will have height $\b$, as desired. 

One can also make the height of the tower go down by activating a sort of wall which prevents the normalizer tower proceeding through:
$$\eqalign{\tri\sp\tri\sp\boxtwotri\sp&\boxfourtri\sp\cdots\cdots\cr
                            &\cboxit{\sp\cboxit{$\sp\square\sp\square\sp$}\sp \cboxit{$\sp\square\sp\square\sp$}\sp}\cr
                            &\cboxit{\sp\cboxit{$\sp\square\sp\square\sp$}\sp \cboxit{$\sp\square\sp\square\sp$}\sp}\cr}$$
By forcing $\square\iso\tri$, it is easy to see that the normalizer tower will not continue past the wall.

In summary, the height of the normalizer tower can be completely controlled by the forcing which makes certain trees isomorphic while retaining their rigidity. The height of the corresponding automorphism tower, therefore, can also be completely controlled, and the proof of the theorem is complete. 

Let me conclude with a discussion of the state of the art with respect to the heights of automorphism towers. In the case of centerless groups, we have Simon Thomas' theorem that the automorphism tower of a group of size $\k$ terminates in fewer than $(2^\k)\plus$ many steps. Of course, this bound cannot be the best possible bound, because there are only $2^\k$ many groups of size $\k$, and so the actual supremum of the towers is strictly below $(2^\k)\plus$. Furthermore, by forcing it is easy to make $2^\k$ as large as desired, and so the ``bound" can be pushed higher and higher. Nevertheless, Thomas' theorem is optimal in the sense that no better upper bound will ever be found:

\theorem.{(Just, Thomas, Shelah 1997)} Suppose that the \GCH\ holds, that $\k$ is a regular uncountable cardinal, that $\cof(\l)>\k$ and that $\a<\l\plus$. Then there is a cardinal-preserving forcing extension in which $2^\k=\l$ and there is a group $G$ of size $\k$ whose automorphism tower has height $\a$. 

\noindent Missing from this analysis are the countable groups. How tall is the automorphism tower of a centerless countable group? 

In the arbitrary case, where non-centerless groups are included, the best uniform upper bound on the height of the automorphism tower of a group is\dots\dots the next inaccessible cardinal. Actually, it is easy to improve slightly on this, by taking the next cardinal $\d$ above the size of the group such that $V_\d$ is a model of \ZFC\ or at least a sufficiently powerful set theory to the prove the Main Theorem above. One has the sense, though, that this is not the right answer. 

Such is our pitiful knowledge even in the case of finite groups! The most we know about an upper bound for the height of the automorphism tower of a finite group is something like the next inaccessible cardinal. Contrast this with the fact that the tallest towers we know of for finite groups have height $\w+n$ for finite $n$. So the true answer lies somewhere between $\w+n$ and the least inaccessible cardinal. How tall is the automorphism tower of a finite group? 

The following questions, to my knowledge, remain open: 

\question. Is there a countable group with an uncountable automorphism tower?

\question. Is there a finite group with an uncountable automorphism tower?

\question. Is there a finite group $G$ such that $G_\w$ is infinite? 

\question. For which ordinals $\g$ is there a group whose tower becomes centerless in exactly
$\g$ many steps? 

\question. Is there a group $G$ whose automorophism tower has height $(2^{\card{G}})\plus$ or more?

And finally, let me state Scott's question, unanswered for 35 years:

\question.{(Scott 1964)} Is the finite part of the automorphism tower of every finite group eventually periodic?

Perhaps Scott's question is peculiar, because no nontrivial instances of this periodicity phenomenon are known. 

\bigskip
\nopagenumbers
\parindent=0pt
\newbox\Article
\newbox\Journal
\newbox\Author
\newbox\Vol
\newbox\No
\newbox\Year
\newbox\Page
\newbox\Book
\newbox\Publisher
\newbox\Pubaddr
\newbox\Key
\newbox\Editor
\newbox\Comment
\newbox\Note
\def\entry#1#2\par{\item{#1\quad}\hskip-1.1em#2\par}
\def\article#1{\setbox\Article=\hbox{\sl #1, }}
\def\journal#1{\setbox\Journal=\hbox{\rm #1 }}
\def\author#1{\setbox\Author=\hbox{\sc #1, }}
\def\vol#1{\setbox\Vol=\hbox{\bf #1 }}
\def\no#1{\setbox\No=\hbox{no. #1 }}
\def\year#1{\setbox\Year=\hbox{\rm({\oldstyle #1}) }}
\def\page#1{\setbox\Page=\hbox{\rm p. #1 }}
\def\book#1{\setbox\Book=\hbox{\it #1, }}
\def\publisher#1{\setbox\Publisher=\hbox{\rm #1, }}
\def\pubaddr#1{\setbox\Pubaddr=\hbox{\rm #1, }}
\def\key#1{\setbox\Key=\hbox{#1}}
\def\editor#1{\setbox\Editor=\hbox{\rm(#1, Ed.) }}
\def\comment#1{\setbox\Comment=\hbox{\rm #1}}
\def\note#1{\setbox\Note=\hbox{\rm #1 }}
\def\ref#1\par{\smallskip{#1
  \entry{\ifhbox\Key\unhbox\Key\else[\ ]\fi}%
  \unhbox\Author\unhbox\Note
  \ifhbox\Book \unhbox\Book\unhbox\Publisher\unhbox\Pubaddr
               \unhbox\Editor\unhbox\Page\unhbox\Year\unhbox\Comment
  \else \unhbox\Article\unhbox\Journal\unhbox\Vol\unhbox\No\unhbox\Editor
        \unhbox\Page\unhbox\Year\unhbox\Comment\fi\par}}

\tenpoint

\ref
\author{Joel David Hamkins and Simon Thomas}
\article{Changing the heights of automorphism towers}
\journal{submitted to the Annals of Pure and Applied Logic}
\key{[1997]}

\ref
\author{Joel David Hamkins}
\article{Every group has a terminating transfinite automorphism tower}
\journal{to appear in the Proceedings of the American Mathematical Society}
\key{[1998]}

\ref
\author{E. Fried and J. Koll\'ar}
\article{Automorphism groups of fields}
\journal{in Universal Algebra (E. T. Schmidt, et al. eds.), Coloq. Math. Soc. Janos Boyali}
\vol{24}
\year{1981}
\page{293-304}
\key{[1981]}

\ref
\author{J. A. Hulse}
\article{Automorphism towers of polycyclic groups}
\journal{Journal of Algebra}
\vol{16}
\year{1970}
\page{347--398}
\key{[1970]}

\ref
\author{Winifried Just, Saharon Shelah and Simon Thomas}
\article{The Automorphism Tower Problem III: Closed Groups of Uncountable Degree}
\comment{Shelah archive \#654}
\key{[1998]}

\ref
\author{Andrew Rae and James E. Roseblade}
\article{Automorphism Towers of Extremal Groups}
\journal{Math. Z.}
\page{70-75}
\year{1970}
\vol{117}
\key{[1970]}

\ref
\author{Scott}
\book{Group Theory}
\year{1964}
\key{[1964]}

\ref
\author{Simon Thomas}
\article{The automorphism tower problem}
\journal{Proceedings of the American Mathematical Society}
\vol{95}
\year{1985}
\page{166--168}
\key{[1985]}

\ref
\author{Simon Thomas}
\article{The automorphism tower problem II}
\journal{Israel J. Math.}
\vol{103}
\year{1998}
\page{93-109}
\key{[1998]}

\ref
\author{H. Wielandt}
\article{Eine Verallgemeinerung der invarianten Untergruppen}
\journal{Math. Z.}
\vol{45}
\year{1939}
\page{209--244}
\key{[1939]}

\bye